\documentclass[10pt]{article} 
\usepackage{amsmath} 
\usepackage{amssymb} 
\usepackage{graphicx} 

\usepackage[numbers]{natbib}

\begin{document}  

\setcounter{figure}{0} 
\setcounter{table}{0} 
\setcounter{footnote}{0} 
\setcounter{equation}{0}  

\vspace*{0.5cm}  

\noindent {\Large \uppercase{Short-axis-mode rotation in complex variables}\footnote{To appear in the Proceedings of the Journ\'ees des Syst\`emes de R\'ef\'erence et de la Rotation Terrestre, 25--27 September 2017, Alicante, Spain}}
\vspace*{0.7cm}  

\noindent\hspace*{1.5cm} M. LARA$^{1,2}$ 
\\
\noindent\hspace*{1.5cm} $^1$ Space Dynamics Group, Technical University of Madrid --- UPM \\
\noindent\hspace*{1.5cm} Pz. Cardenal Cisneros 3, 28040 Madrid, Spain \\
\noindent\hspace*{1.5cm} $^2$ GRUCACI, University of La Rioja \\
\noindent\hspace*{1.5cm} C/ Madre de Dios, 53, 26006 Logro\~no, Spain \\
\noindent\hspace*{1.5cm} e-mail: mlara0@gmail.com  \\

\vspace*{0.5cm}

\noindent {\large ABSTRACT.} Decomposition of the free rigid body Hamiltonian into a ``main problem'' and a perturbation term provides an efficient integration scheme that avoids the use of elliptic functions and integrals. In the case of short-axis-mode rotation, it is shown that the use of complex variables converts the integration of the torque-free motion by perturbations into a simple exercise of polynomial algebra that can also accommodate the gravity-gradient perturbation when the rigid body rotation is close enough to the axis of maximum inertia.
\vspace*{1cm}  

\noindent {\large 1. \uppercase{Introduction}}  
\smallskip  

The rotation of a rigid body in the absence of external torques is known to be integrable \cite{Whittaker1917,Golubev1960}. However, because the solution depends on elliptical integrals and elliptic functions, in practical applications the closed form solution is customarily replaced by corresponding expansions in trigonometric functions truncated to some order. In particular, useful expansions by Kinoshita apply when either the triaxiality of the rigid body is small \cite{Kinoshita1972} or the rotation is close to either the axis of maximum or minimum momentum of inertia \cite{Kinoshita1992} ---the order of these expansions was later extended by other authors \cite{SouchayFolgueiraBouquillon2003}. Alternatively, the expansions in trigonometric functions of the solution of the free rigid body can be directly constructed using perturbation theory \cite{FerrerLara2010,Lara2014}, an approach that systematizes the computation of higher orders of the expansions and eases the construction of perturbation solutions in the presence of external torques \cite{LaraFukushimaFerrer2011}.
\par

Perturbation approaches to the torque-free motion of a rigid body start from the decomposition of the free rigid body Hamiltonian in Andoyer variables into a ``main problem'' and a perturbation term. When the triaxiality is small, the selection of either the axisymmetric case or the spherical rotor as the main problem results in a zeroth order Hamiltonian that only depends on momenta of the canonical set of Andoyer variables, which, therefore, are directly action-angle variables, a fact that simplifies application of the perturbation method \cite{FerrerLara2010}. On the contrary, in the case of short-axis-mode (SAM) rotation, the main problem involves the same Andoyer variables as the free rigid body problem, which include coordinates as well as their conjugate momenta, and, in consequence, are not action-angle variables of the zeroth order Hamiltonian. Therefore, in order to set up an efficient perturbative integration scheme, a preliminary complete reduction of the main problem of SAM rotation in action-angle variables, which is achieved solving the Hamilton-Jacobi equation, is required \cite{Lara2014}.  
\par

On the other hand, it will be shown that the perturbative arrangement of the free rigid body Hamiltonian in the case of SAM rotation is immediately disclosed when using non-singular variables of Poincar\'e type \cite{HenrardMoons1978}. In these variables, the free rigid body Hamiltonian takes the form of the Hamiltonian of the simple harmonic oscillator disturbed by additional quartic polynomial terms. A following transformation to complex variables converts the integration of the free rigid body Hamiltonian by perturbations into a simple exercise of polynomial algebra. The polynomial structure of the perturbation is not preserved, in general, when perturbation torques are taken into account. However, when the rotation is close enough to the axis of maxima inertia, it is shown that the gravity-gradient perturbation can be easily tackled within the same perturbative scheme.
\vspace*{0.7cm}  

\newpage
\noindent {\large 2. \uppercase{Perturbative arrangement}}  
\smallskip  

Andoyer's \cite{Andoyer1923} arrangement of the rigid body Hamiltonian is
\begin{equation}\label{hacheA}
\mathcal{H}_0=\frac{M^2}{2C}\left[1+\alpha\left(1-\frac{N^2}{M^2}\right)-\alpha\,\beta\left(1-\frac{N^2}{M^2}\right)\cos2\nu\right],
\end{equation}
where $(\lambda,\mu,\nu,\Lambda,M,N)$ are the usual Andoyer variables defining the node of the invariable plane on the inertial $x$, $y$ plane, the node of the equatorial plane of the rigid body on the invariable plane, the component of the rotation of the body around its axis of maximum inertia, the projection of the total angular momentum along the inertial $z$ axis, the total angular momentum, and its projection along the body axis of maximum inertia, respectively. The relations
\begin{equation} \label{alfabeta}
\alpha\,(1+\beta)=\frac{C}{A}-1,
\qquad
\alpha\,(1-\beta)=\frac{C}{B}-1,
\end{equation}
define the physical parameters $\alpha$ and $\beta$ as a function of the principal moments of inertia $A\le{B}\le{C}$. When the \emph{triaxiality coefficient} $\beta$ is small, Eq.~(\ref{hacheA}) admits a perturbative arrangement in which the zeroth order term corresponds to an axisymmetric body, whose Hamiltonian is completely reduced, and the perturbation, which is due to the triaxiality, depends on the angle $\nu$. This perturbative arrangement eases the computation of a perturbation solution in trigonometric functions that matches Kinoshita's \cite{Kinoshita1972} series expansion of the closed form solution in powers of $\beta$ \cite{FerrerLara2010}, and can be easily extended to any order of $\beta$.
\par


An alternative perturbative arrangement has been recently proposed for rigid bodies rotating close to its axis of maximum inertia, irrespective of its triaxiality \cite{Lara2014}. In that case $N\approx{M}$ and, therefore, $\frac{1}{2}(1-N/M)=\sin^2\frac{1}{2}J\ll1$, where $J$ is the inclination angle between the invariable plane and the equatorial plane of the rigid body, a fact that motivates reorganization of Eq.~(\ref{hacheA}) in the form
\begin{equation}\label{hacheP}
\mathcal{H}_0=\mathcal{M}+\varepsilon\,\mathcal{P},
\end{equation}
where 
\begin{equation}\label{main}
\mathcal{M}=\frac{M^2}{2C}\left[1+2\alpha\,\left(1-\frac{N}{M}\right)\left(1-\beta\cos2\nu\right)\right],
\end{equation}
is taken as the integrable part, $\varepsilon$ is a formal small parameter, and
\begin{equation}\label{disturbing}
\mathcal{P}=-\frac{M^2}{2C}\,\alpha\left(1-\frac{N}{M}\right)^2\left(1-\beta\cos2\nu\right),
\end{equation}
is a perturbation $|\mathcal{P}|\ll\mathcal{M}$.
\par

Now, the zeroth order Hamiltonian (\ref{main}), which has been dubbed as the main problem of SAM rotation, involves the same variables as the free rigid body Hamiltonian (\ref{hacheA}). For a perturbation approach, the Hamiltonian is customarily reformulated in new action-angle variables such that the zeroth-order term is completely reduced to a function of only the new momenta. The complete reduction of Eq.~(\ref{main}), which was achieved in \cite{Lara2014} by solving the Hamilton-Jacobi equation in the style of \cite{FerrerLara2010b}, becomes trivial when using non-singular variables.
\par


Andoyer variables are singular for $N=M$, a case in which $\nu$ is not defined. However, this singularity is virtual \citep{Henrard1974}, and is easily avoided using non-singular variables of the Poincar\'e type \cite{HenrardMoons1978}. 
Thus, the canonical transformation 
\begin{eqnarray} \label{xn}
x &=& -\sqrt{2(M-N)}\sin\nu, \\ \label{Xn}
X &=& \sqrt{2(M-N)}\cos\nu, \\ \label{yn}
y &=& \mu+\nu, \\ \label{Yn}
Y &=& M,
\end{eqnarray}
converts Eq.~(\ref{main}) into
\begin{equation} \label{HamHO}
\mathcal{M} = \frac{Y^2}{2C} + \frac{Y}{C}\alpha(1-\beta)
\frac{1}{2}\left(X^2+\omega^2x^2\right),
\end{equation}
with
\begin{equation} \label{wb}
\omega=\sqrt\frac{1+\beta}{1-\beta}, \qquad \beta=\frac{\omega^2-1}{\omega^2+1}.
\end{equation}
\par

One easily recognizes in Eq.~(\ref{HamHO}) the Hamiltonian of a harmonic oscillator of (non-dimensional) frequency $\omega$, and it is well known that the Poincar\'e transformation $(\ell,L,\omega)\longrightarrow(x,X)$ given by
\begin{equation} \label{itns}
x=\sqrt{2L/\omega}\sin\ell, \qquad
X=\sqrt{2\omega\,L}\cos\ell,
\end{equation}
completely reduces this Hamiltonian to a function of only the momentum $L$. In this way ---formulation of the main problem Hamiltonian in the nonsingular variables in Eqs.~(\ref{xn})--(\ref{Yn}) followed by the Poincar\'e transformation in Eq.~(\ref{itns})--- the computation of the action-angle variables of the main problem of the SAM rotation carried out in \citep{Lara2014} is dramatically abridged.
\par

On the other hand, the use of action-angle variables, while customary, is not a requirement in perturbation theory. Indeed, in view of Eq.~(\ref{disturbing}) takes the form of a quartic polynomial in the nonsingular variables of Eq.~(\ref{xn})--(\ref{Yn}), viz.
\begin{equation} \label{disturbingquartic}
\mathcal{P} =-\frac{\alpha}{8C}\left[(1+\beta)x^4+2x^2X^2+(1-\beta)X^4\right],
\end{equation}
the perturbation solution can be directly constructed in Cartesian variables. Moreover, it is known that the use of complex variables makes the procedure very efficient \citep{Kummer1976,GiorgilliGalgani1978}.
\vspace*{0.7cm}  

\noindent {\large 3. \uppercase{Perturbation solution in complex variables}}  
\smallskip  

If the transformation
\begin{equation} \label{tocomplex}
x=\frac{1}{\sqrt{2\omega}}(u-i U), \quad
X=\sqrt{\frac{\omega}{2}}(U-iu), \qquad
y=\sqrt{\gamma}\,v, \quad
Y=\frac{1}{\sqrt\gamma}V,
\end{equation}
where $i=\sqrt{-1}$ and
\[
\gamma=\alpha\sqrt{1-\beta^2}=\sqrt{\left(\frac{C}{A}-1\right)\left(\frac{C}{B}-1\right)},
\] 
which is canonical, is now performed, Eq.~(\ref{HamHO}) is rewritten in the real $(v,V)$ and complex $(u,U)$ variables in the form
\begin{equation} \label{MainComplex}
\mathcal{M}=\frac{V^2}{2C\gamma}-\frac{V\sqrt{\gamma}}{C}iuU,
\end{equation}
whereas Eq.~(\ref{disturbingquartic}) takes the form
\begin{equation} \label{Pcomplex}
\mathcal{P}=\frac{\alpha}{4C}\left[2 u^2 U^2-i \beta  \left(u^3 U-u U^3\right)\right].
\end{equation}

The Lie derivative $\mathcal{L}_\mathcal{M}$ associated to Eq.~(\ref{MainComplex}), is given by the Poisson bracket operator $\mathcal{L}_\mathcal{M}=\{\quad;\mathcal{M}\}$, viz.
\begin{equation} \label{Liefull}
\mathcal{L}_\mathcal{M}=
\frac{\sqrt\gamma}{C}\left[Vi\left(U\frac{\partial}{\partial{U}}-u\frac{\partial}{\partial{u}}\right)
+\left(\frac{V}{\gamma^{3/2}}+iuU\right)\frac{\partial}{\partial{v}}\right],
\end{equation}
and the partial differential equation $\mathcal{L}_\mathcal{M}(W_n)=\widetilde{M}_{n}-{M}_{0,n}$ must be solved at each order $n$ of the perturbation theory to compute the corresponding term $W_n$ of the generating function. Terms $\widetilde{M}_{n}$ are known from previous computations whereas terms ${M}_{0,n}$ are chosen to pertain to the kernel of the Lie derivative, viz. $\mathcal{L}_\mathcal{M}(W_n)=0$ (see \cite{Deprit1969} for details).
\par

However, because Eq.~(\ref{Pcomplex}) does not depend on $v$, when dealing with the torque free motion as a perturbation problem one can assume that the generating function is independent of $v$. Hence,
\begin{equation} \label{Liefree}
\mathcal{L}_\mathcal{M}=
\frac{\sqrt\gamma}{C}Vi\left(U\frac{\partial}{\partial{U}}-u\frac{\partial}{\partial{u}}\right).
\end{equation}
Then, for any integers $j\ge0$ and $k\ge0$,
\[
\mathcal{L}_\mathcal{M}(u^j\,U^k) = \frac{\sqrt\gamma}{C}Vi(k-j)u^jU^k,
\]
and, therefore, $\mathcal{L}_\mathcal{M}(u^j\,U^k)=0$ requires that $j=k$. That is,  the kernel of the Lie operator is composed of monomials of the form $(u\,U)^j$, whereas all other monomials $u^j\,U^k$, $j\ne{k}$, pertain to the image. Therefore, the solution of the homological equation becomes trivial in complex variables. Indeed, any monomial $q_{j,k}u^j\,U^k$, $j\ne{k}$, where $q_{j,k}$ is a numeric coefficient, contributes a term
\[
i\frac{C}{\sqrt\gamma}\frac{q_{j,k}}{(j-k)V}u^jU^k,
\]
to the generating function.


The procedure starts from writing the free rigid body Hamiltonian like the Taylor series expansion
\[
\mathcal{H}_0=\sum_{n\ge0}\frac{\varepsilon^n}{n!}H_{n,0}(u,U),
\]
where $H_{0,0}\equiv\mathcal{M}$ as given in Eq.~(\ref{MainComplex}), $H_{1,0}\equiv\mathcal{P}$ as given by Eq.~(\ref{Pcomplex}), and $H_{n,0}=0$ for $n\ge2$. Straightforward computations lead to the normalized Hamiltonian, in new, prime variables
\[
\mathcal{K}=\frac{\sqrt{\gamma}}{2C}V^2\left[\frac{1}{\gamma^{3/2}}+
\sum_{n\ge0}\frac{\alpha^n}{\gamma^{n/2}}p_n\left(\frac{u'U'}{V}\right)^{n+1}\right],
\]
where $p_n$ are polynomials in the triaxiality coefficient $\beta$. The first few triaxiality polynomials are $p_0=-2i$, $p_1=1$, $p_2=\frac{1}{2}\beta^2$, $p_3=\frac{5}{8}\beta^2$, and
\begin{eqnarray*}
p_4 &=& \frac{3}{32}\beta^2\left(3\beta^2+8\right) \\
p_5 &=& \frac{7}{32}\beta^2\left(5\beta^2+4\right) \\
p_6 &=& \frac{1}{128}\beta^2\left(45\beta^4+354\beta^2+128\right) \\
p_7 &=& \frac{9}{1024}\beta^2\left(265\beta^4+650\beta^2+128\right) \\
p_8 &=& \frac{5}{8192}\beta^2\left(953\beta^6+14888\beta^4+17120\beta^2+2048\right) \\
p_9 &=& \frac{11}{8192}\beta^2\left(4075\beta^6+20212\beta^4+13104\beta^2+1024\right)
\end{eqnarray*}
which, as expected, are the same as those in Table 2 of \citep{Lara2014} after adjusting subindices and scaling by $\beta^2$.
\par

The transformation from prime to original variables
\[
u=u'+\sum_{n\ge1}\frac{\varepsilon^n}{n!}u_{0,n}(u',U'), \qquad
U=U'+\sum_{n\ge1}\frac{\varepsilon^n}{n!}U_{0,n}(u',U'),
\]
is obtained by successive evaluations of Deprit's triangle
\begin{equation} \label{triangle}
f_{n,q}=f_{n+1,q-1}+
\sum_{0\le{m}\le{n}}{n\choose{m}}\,\{f_{n-m,q-1};S_{m+1}\},
\end{equation}
using the generating function $\mathcal{S}=\sum_{m\ge0}(\varepsilon^m/m!)\,S_{m+1}$ where
\[
S_m=\frac{(m-1)!\alpha^m}{\gamma^{m/2}V^m}\beta\left(u^2+U^2\right)s_m
\]
and the first few $s_m$ are
\begin{eqnarray*}
s_1 &=& \frac{i}{8}uU, \\
s_2 &=& \frac{i}{4}u^2U^2, \\
s_3 &=& \frac{1}{64}\left[{24i\left(\beta ^2+2\right)uU+5\beta(U^2-u^2)}\right]u^2U^2, \\
s_4 &=& \frac{3}{64}\left[{2i\left(57\beta^2+32\right) u U+\beta\left(9\beta^2+20\right)(U^2-u^2)}\right]u^3U^3, \\
s_5 &=& \frac{i}{64} \left[
2 \left(343\beta^4+2024\beta^2+480\right) u^2U^2
-11 \beta^2\left(\beta^2+2\right)\left(u^4+U^4\right) \right. \\
&& \left.-6i\beta\left(147\beta^2+100\right)uU\left(U^2-u^2\right)
\right]u^3U^3.
\end{eqnarray*}

\newpage
\noindent {\large 4. \uppercase{Gravity gradient}}  
\smallskip  

The perturbation approach based on the main problem of SAM rotation is also feasible for motion under external torques. In the particular case of gravity-gradient perturbations due to a distant body, the majority of perturbation terms are factored by $\sin{J}$ or $\sin^2J$ (see \cite{LaraFukushimaFerrer2011}, for instance). Besides, typical values of the gravity-gradient perturbation for solar system bodies are small when compared with the torque-free rotation, say below $10^{-6}$. Therefore, in those cases in which the inclination angle $J$ is small, terms $\mathcal{O}(\sin\frac{1}{2}J)$ can be neglected.
\par

Then, if one makes the simplifying assumption that the disturbing body moves with Keplerian motion, and takes its orbital plane as the inertial plane, the only relevant terms of the gravity-gradient perturbation, in the new variables, are simply
\begin{eqnarray} \label{gragra}
\mathcal{D}&=&-\frac{n^2}{4}\,\frac{a^3}{r^3}\bigg\{\!
\left(C-\frac{A+B}{2}\right)\left(2-3s^2+3s^2\cos2\vartheta\right)
\\ \nonumber
&& +\frac{3}{4}(B-A)\left[(1-c)^2\cos(2y-2\vartheta)+2s^2\cos2y+(1+c)^2\cos(2y+2\vartheta)\right]\!\bigg\},
\end{eqnarray}
where $\vartheta=\lambda-\theta$, $r$ and $\theta$ are polar coordinates, $a$ is orbit semimajor axis, $n$ is orbital mean motion, $c\equiv\cos{i}$, $s\equiv\sin{i}$, and $i=\arccos(\Lambda/M)$ is the inclination angle between the orbital plane and the invariable plane.
\par

Now, the full Lie derivative in Eq.~(\ref{Liefull}) is involved in the solution of the homological equation.
Note that the first summand in the square brackets of Eq.~(\ref{Liefull}) vanishes for terms of the form $F(v,V,uU)$. Hence, because Eq.~(\ref{gragra}) is made of terms of this type, and in view of the form of the second summand in the square brackets of Eq.~(\ref{Liefull}), which only includes a factor $uU$, the homological equation is easily solved. Indeed, if we choose the new Hamiltonian term
\[
\langle\mathcal{D}\rangle=\frac{1}{2\pi}\int_0^{2\pi}\mathcal{D}\,\mathrm{d}y=
-\frac{n^2}{4}\,\frac{a^3}{r^3}\left(C-\frac{A+B}{2}\right)\left(2-3s^2+3s^2\cos2\vartheta\right),
\]
and assume that there is no coupling with the previous terms of the perturbation theory, a particular solution of the homological equation for the order $n$ corresponding to this term, is
\begin{eqnarray*}
S_n &=& -\frac{n^2}{4}\,\frac{a^3}{r^3}\frac{3}{4}(B-A)\frac{3}{16}\frac{C\sqrt{\gamma}}{i\gamma^{3/2}uU-V}\times \\
&&\left[(1-c)^2\sin(2y-2\vartheta)+2s^2\sin2y+(1+c)^2\sin(2y+2\vartheta)\right].
\end{eqnarray*}
\vspace*{0.7cm}  

\noindent {\large 5. \uppercase{Conclusions}}  
\smallskip  

Short-axis mode rotation of a free rigid body is naturally decomposed into a main problem and a perturbation, a fact that leads to the straightforward integration of the rotation by perturbation series. When using non-singular variables of the Poincar\'e type, the main problem has the form of a harmonic oscillator, whose frequency is related to the triaxiality of the rigid body, whereas the perturbation is a quartic polynomial. Then, the use of complex variables makes the construction of the perturbation solution trivial. The polynomial character of the perturbation does not persist, in general, when the motion is affected by external torques. However, when the rotation is close to the axis of maximum inertia, the gravity-gradient perturbation can also be approached in complex variables.
\vspace*{0.7cm}  

\noindent {\uppercase{Acknowledgements}}  
\smallskip  

Partial support by the Spanish State Research Agency and the European Regional Development Fund under Projects ESP2016-76585-R and ESP2013-41634-P (AEI/ER\-DF, EU) are recognized. Discussions with A.~Escapa, University of Le\'on, motivated the development of this work.


\vspace*{0.7cm}  

\newpage

\noindent {\large 7. REFERENCES} 
{  

\leftskip=5mm 
\parindent=-5mm  
\smallskip  

\begin{list}{}{}

\bibitem[{{Andoyer}(1923)}]{Andoyer1923}
{Andoyer} MH (1923) {Cours de m\'ecanique c\'eleste. Tome I.} {Paris:
  Gauthier-Villars, 438 S. (1923).}

\bibitem[{{Deprit}(1969)}]{Deprit1969}
{Deprit} A (1969) Canonical transformations depending on a small parameter.
  Celestial Mechanics 1(1):12--30 

\bibitem[{{Ferrer} and {Lara}(2010{\natexlab{a}})}]{FerrerLara2010b}
{Ferrer} S, {Lara} M (2010{\natexlab{a}}) {Families of Canonical
  Transformations by Hamilton-Jacobi-Poincar\'e Equation. Application to
  Rotational and Orbital Motion}. Journal of Geometric Mechanics 2(3):223--241,

\bibitem[{{Ferrer} and {Lara}(2010{\natexlab{b}})}]{FerrerLara2010}
{Ferrer} S, {Lara} M (2010{\natexlab{b}}) {Integration of the Rotation of an
  Earth-like Body as a Perturbed Spherical Rotor}. The Astronomical Journal
  139(5):1899--1908 

\bibitem[{{Giorgilli} and {Galgani}(1978)}]{GiorgilliGalgani1978}
{Giorgilli} A, {Galgani} L (1978) {Formal integrals for an autonomous
  Hamiltonian system near an equilibrium point}. Celestial Mechanics
  17:267--280 

\bibitem[{Golubev(1960)}]{Golubev1960}
Golubev V (1960) Lectures on Integration of the Equations of Motion of a Rigid
  Body about a Fixed Point. Israel Program for Scientific Translations,
  S.~Monson, Jerusalem

\bibitem[{{Henrard}(1974)}]{Henrard1974}
{Henrard} J (1974) Virtual singularities in the artificial satellite theory.
  Celestial Mechanics 10(4):437--449 

\bibitem[{Henrard and Moons(1978)}]{HenrardMoons1978}
Henrard J, Moons M (1978) Hamiltonian theory of the libration of the moon. In:
  Szebehely VG (ed) Dynamics of planets and satellites and theories of their
  motion, Proceedings of the International Astronomical Union colloquium
  no.~41, D.~Reidel Publishing Company, Dordrecht: Holland / Boston: U.S.A.,
  Astrophysics and Space Science Library, vol~72, pp 125--135

\bibitem[{{Kinoshita}(1972)}]{Kinoshita1972}
{Kinoshita} H (1972) {First-Order Perturbations of the Two Finite Body
  Problem}. Publications of the Astronomical Society of Japan 24:423--457

\bibitem[{{Kinoshita}(1992)}]{Kinoshita1992}
{Kinoshita} H (1992) {Analytical expansions of torque-free motions for short
  and long axis modes}. Celestial Mechanics and Dynamical Astronomy
  53(4):365--375 

\bibitem[{{Kummer}(1976)}]{Kummer1976}
{Kummer} M (1976) {On resonant non linearly coupled oscillators with two equal
  frequencies}. Communications in Mathematical Physics 48:53--79,

\bibitem[{Lara(2014)}]{Lara2014}
Lara M (2014) Short-axis-mode rotation of a free rigid body by perturbation
  series. Celestial Mechanics and Dynamical Astronomy 118(3):221--234,

\bibitem[{{Lara} et~al(2011){Lara}, {Fukushima}, and
  {Ferrer}}]{LaraFukushimaFerrer2011}
{Lara} M, {Fukushima} T, {Ferrer} S (2011) {Ceres' rotation solution under the
  gravitational torque of the Sun}. Monthly Notices of the Royal Astronomical
  Society 415(1):461--469 

\bibitem[{{Souchay} et~al(2003){Souchay}, {Folgueira}, and
  {Bouquillon}}]{SouchayFolgueiraBouquillon2003}
{Souchay} J, {Folgueira} M, {Bouquillon} S (2003) {Effects of the Triaxiality
  on the Rotation of Celestial Bodies: Application to the Earth, Mars and
  Eros}. Earth Moon and Planets 93(2):107--144,

\bibitem[{{Whittaker}(1917)}]{Whittaker1917}
{Whittaker} ET (1917) {A Treatise on the Analytical Dynamics of Particles and
  Rigid Bodies}, 2nd edn. Cambridge University Press

\end{list}
%
%
%
%

} 

\end{document}